\documentclass[12pt,a4paper,twoside]{article}

\usepackage{amsfonts,amsmath,amssymb,amsthm,latexsym}

\usepackage{mathptmx}

\usepackage[latin1]{inputenc} 
\usepackage{eucal}            
\usepackage{bbm}              
\usepackage[ps]{xypic}        
\usepackage{diagxy}           

\newtheorem{proposition}{Proposition}[section]

\newtheorem{theorem}[proposition]{Theorem}

\theoremstyle{definition}
\newtheorem{definition}[proposition]{Definition}
\newtheorem{example}[proposition]{Example}

\theoremstyle{remark}
\newtheorem{remark}[proposition]{Remark}

\renewcommand{\mathbb}[1]{\mathbbm{#1}} 

\newcommand{\I}          {\mathrm{i}}

\newcommand{\cc}[1]      {\overline{{#1}}}              
\newcommand{\id}         {\operatorname{\mathsf{id}}}   
\newcommand{\Hom}        {\operatorname{\mathsf{Hom}}}   
   
\newcommand{\SP}[1]      {\left\langle{#1}\right\rangle} 
\newcommand{\ring}[1]    {\mathsf{#1}}                 
\newcommand{\Unit}       {\mathbb{1}}                  

\newcommand{\acts}       {\mathbin{\triangleright}}
\newcommand{\Bimod}[5] {\sideset{^{\scriptscriptstyle{#1}}_{\scriptscriptstyle{#2}}}{^{\scriptscriptstyle{#4}}_{\scriptscriptstyle{#5}}}{\operatorname{#3}}}
\newcommand{\EA}   {\Bimod{}{}{\mathcal{E}}{}{\mathcal{A}}}
\newcommand{\EpA}  {\Bimod{}{}{\mathcal{E}}{\prime}{\mathcal{A}}}
\newcommand{\BEA}  {\Bimod{}{\mathcal{B}}{\mathcal{E}}{}{\mathcal{A}}}
\newcommand{\CFB}  {\Bimod{}{\mathcal{C}}{\mathcal{F}}{}{\mathcal{B}}}
\newcommand{\AAA}  {\Bimod{}{\mathcal{A}}{\mathcal{A}}{}{\mathcal{A}}}
\newcommand{\sweedler}[1] {{\scriptscriptstyle{(#1)}}}
\newcommand{\twist}[1]   {\textbf{#1}}

\newcommand{\IP}[4]{{\,}_{\scriptscriptstyle{#2}\!\!}\left\langle{{#1}}\right\rangle^{\scriptscriptstyle{#3}}_{\scriptscriptstyle{#4}}}

\newcommand{\SPA}[1]     {\IP{{#1}}{}{}{\mathcal{A}}}
\newcommand{\SPEA}[1]    {\IP{{#1}}{}{\mathcal{E}}{\mathcal{A}}}
\newcommand{\SPEpA}[1]    {\IP{{#1}}{}{\mathcal{E}'}{\mathcal{A}}}
\newcommand{\BSPE}[1]    {\IP{{#1}}{\mathcal{B}}{\mathcal{E}}{}}
\newcommand{\SPFB}[1]    {\IP{{#1}}{}{\mathcal{F}}{\mathcal{B}}}
\newcommand{\SPFEA}[1]   {\IP{{#1}}{}{\mathcal{F}\otimes\mathcal{E}}{\mathcal{A}}}

\newcommand{\tensor}[1][{}]{\mathbin{\otimes_{\scriptscriptstyle{#1}}}}
\newcommand{\tensM}[1][{}] {\mathbin{\widehat{\otimes}_{\scriptscriptstyle{#1}}}}

\newcommand{\Rep}[1][{}]  {\sideset{^*}{_{#1}}{\operatorname{\textrm{-}\mathsf{Rep}}}}
\newcommand{\sMod}[1][{}] {\sideset{^*}{_{#1}}{\operatorname{\textrm{-}\mathsf{Mod}}}}

\newcommand{\BiMods}      {\sideset{}{^*_{}}{\operatorname{\mathsf{Bimod}}}}
\newcommand{\BiModstr}    {\sideset{}{^{\mathsf{str}}}{\operatorname{\mathsf{Bimod}}}}
\newcommand{\BiModsH}     {\sideset{}{^*_{H}}{\operatorname{\mathsf{Bimod}}}}
\newcommand{\BiModstrH}   {\sideset{}{^{\mathsf{str}}_{H}}{\operatorname{\mathsf{Bimod}}}}

\newcommand{\Pic}      {\operatorname{\mathsf{Pic}}}
\newcommand{\PicH}      {\sideset{}{_{H}}{\operatorname{\mathsf{Pic}}}}
\newcommand{\StrPic}   {\sideset{}{^{\mathrm{str}}}{\operatorname{\mathsf{Pic}}}}
\newcommand{\starPic}  {\sideset{}{^*}{\operatorname{\mathsf{Pic}}}}
\newcommand{\StrPicH}  {\sideset{}{^{\mathrm{str}}_{H}}{\operatorname{\mathsf{Pic}}}}
\newcommand{\starPicH} {\sideset{}{^*_{H}}{\operatorname{\mathsf{Pic}}}}

\renewcommand{\thefootnote}{\fnsymbol{footnote}}

\title{\MakeUppercase{Covariant Strong Morita Theory of Star Product
    Algebras}\thanks{Contributions to the proceedings of the NoMap
    conference, 2008.}}
\author{Stefan Waldmann}
\date{}
\renewcommand{\date}{\vspace{-5mm}}
\begin{document}
\maketitle \vspace*{-3mm}\relax
\renewcommand{\thefootnote}{\arabic{footnote}}

\noindent 
\textit{\small Fakultät für Mathematik und Physik, Hermann Herder
  Straße 3, D 79104 Freiburg, Germany\\
  email: Stefan.Waldmann@physik.uni-freiburg.de}

%
%

\begin{abstract}
\noindent
In this note we recall some recent progress in understanding the
representation theory of $^*$-algebras over rings $\ring{C} =
\ring{R}(\I)$ where $\ring{R}$ is ordered and $\I^2 = -1$. The
representation spaces are modules over auxiliary $^*$-algebras with
inner products taking values in this auxiliary $^*$-algebra. The ring
ordering allows to implement positivity requirements for the inner
products. Then the representations are required to be compatible with
the inner product. Moreover, one can add the notion of symmetry in
form of Hopf algebra actions. For all these notions of representations
there is a well-established Morita theory which we review. The core of
each version of Morita theory is the corresponding Picard groupoid for
which we give tools to compute and determine both the orbits and the
isotropy groups.
\end{abstract}

%
%

\section{Introduction}
\label{sec:Introduction}

Deformation quantization \cite{bayen.et.al:1978a} has proven to be a
physically reasonable and mathematically rich approach to quantization
of finite-dimensional classical mechanical systems modelled on phase
spaces being symplectic or even general Poisson manifolds. Many
questions like existence and classification of star products have been
successfully answered and non-trivial applications to other branches
of mathematics like index theorems have been found. To name a few
highlights see \cite{kontsevich:2003a, nest.tsygan:1995a,
  fedosov:1996a, dewilde.lecomte:1983b, bertelson.cahen.gutt:1997a}
and confer to \cite{waldmann:2007a} for a gentle introduction with
further and detailed references.

While the structure of the deformed algebras, playing the role of
observables in deformation quantization, is by now very well
understood, this is not yet enough for the original physical purpose:
One also has to develop a physically reasonable notion for states
allowing for an implementation of the superposition principle. The
solution of this problem is to proceed analogously as in the theory of
$C^*$-algebras: states are identified with positive linear
functionals, the superposition principle is then implemented in the
GNS representations of the star product algebra arising from such
positive functionals. This program has been pursued successfully in a
sequence of articles, see in particular \cite{bursztyn.waldmann:2005b,
  waldmann:2005b, waldmann:2002a, bordemann.waldmann:1998a} as well as
\cite[Chap.~7]{waldmann:2007a} for an overview and more references.

The aim of this note is now to review recent developments in studying
the representation theory of star product algebras in the presence of
a \emph{symmetry}. The physical relevance of this question should be
clear.

The results presented in this review are mainly based on the work
\cite{jansen.waldmann:2006a}, see also \cite{jansen:2006a} and the
forthcoming project \cite{jansen.neumaier.waldmann:2005a:pre}. We
start explaining the basic ingredients in Section~\ref{sec:SetUp}.
Then we define the adapted notions of representations in
Section~\ref{sec:RepresentationTheories} and introduce the relevant
tensor products for Morita theory in
Section~\ref{sec:TensorProductsMoritaTheory}. Based on these tensor
products we explain the induced notions of Picard (bi-)groupoids in
Section~\ref{sec:PicardGroupoids} and comment on some of the Morita
invariants arising from actions of the Picard groupoids. Finally, we
conclude this review with some further remarks in
Section~\ref{sec:FurtherRemarks}.

\noindent
\textbf{Acknowledgements:} The author would like to thank the
organizers of the NoMap 2008 for their excellent organization and the
participants for many valuable comments and discussions.

%
%

\section{The set-up}
\label{sec:SetUp}

The framework we discuss here is slightly more general than actually
needed for star product algebras. First we consider an \emph{ordered
  ring} $\ring{R}$ which in the examples will mainly be either
$\mathbb{R}$ or $\mathbb{R}[[\lambda]]$ both with their usual ordering
structure. The positivity in $\ring{R}$ will then be used to induce
notions of positivity for linear functionals, algebra elements, inner
products, etc. Next, we take the ring extension $\ring{C} =
\ring{R}(\I)$ by a square root of $-1$, i.e. $\I^2 = -1$ to have a
replacement for the complex numbers. In the main two examples this
results in $\mathbb{C}$ and $\mathbb{C}[[\lambda]]$, respectively.

The algebras we consider in this note are \emph{$^*$-algebras}
$\mathcal{A}$ over $\ring{C}$, i.e. associative algebras over
$\ring{C}$ equipped with a antilinear involutive antiautomorphism, the
\emph{$^*$-involution}, written as $a \mapsto a^*$ for $a \in
\mathcal{A}$. The notion of \emph{symmetry} we are interested is based
on a Hopf algebra $H$ which we require to be a Hopf $^*$-algebra, i.e.
a Hopf algebra which is a $^*$-algebra such that the comultiplication
$\Delta$ and the counit $\epsilon$ are $^*$-homomorphisms and such
that $S {}^* S {}^* = \id$ for the antipode $S$. In fact, the
requirement $S {}^* S {}^* = \id$ is superfluous and follows already
from the remaining properties of a Hopf $^*$-algebra, see e.g.
\cite{klimyk.schmuedgen:1997a} for more details on Hopf algebras with
$^*$-involutions. Then all $^*$-algebras $\mathcal{A}$ are though to
carry a \emph{left $^*$-action} of the Hopf $^*$-algebra, i.e. a left
$H$-module structure such that
\begin{equation}
    \label{eq:LeftAction}
    g \acts(ab) = (g_\sweedler{1} \acts a) (g_\sweedler{2} \acts b),
    \quad
    g \acts \Unit_{\mathcal{A}} = \epsilon(g) \Unit_{\mathcal{A}},
    \quad
    \textrm{and}
    \quad
    (g \acts a)^* = S(g)^* \acts a^*,
\end{equation}
where $\Delta(g) = g_\sweedler{1} \otimes g_\sweedler{2}$ is the usual
Sweedler notation and $g \in H$ and $a, b \in \mathcal{A}$. The main
examples of interest in deformation quantization are the group
algebras $H = \ring{C}[G]$ with $^*$-involution specified by $g^* =
g^{-1}$ and the usual Hopf algebra structure as well as the
complexified universal enveloping algebras $U(\mathfrak{g})
\otimes_{\ring{R}} \ring{C}$ of a Lie algebra $\mathfrak{g}$ over
$\ring{R}$ with $^*$-involution determined by $\xi^* = -\xi$ for $\xi
\in \mathfrak{g}$ and the usual Hopf algebra structure.

We call such a $^*$-action \emph{inner} if there is a
$^*$-homomorphism, the \emph{momentum map} $J: H \longrightarrow
\mathcal{A}$, such that
\begin{equation}
    \label{eq:MomentumMap}
    g \acts a = J(g_\sweedler{1}) a J(S(g_\sweedler{2})).
\end{equation}
In general, it is a quite strong requirement to have a momentum map
and many interesting cases do \emph{not} allow for such a $^*$-action.
In particular, if the action is non-trivial and $\mathcal{A}$ is
commutative, then there can be no momentum map. However, there is one
important exception in deformation quantization: assume that we have a
real Lie algebra $\mathfrak{g}$ acting on a manifold $M$ by vector
fields. Then we have an action by derivations on $C^\infty(M)$ via the
Lie derivative. Assume now that we have a star product $\star$ which
has a quantum momentum map. This is a linear map $J: \mathfrak{g}
\longrightarrow C^\infty(M)[[\lambda]]$ such that
\begin{equation}
    \label{eq:QuantumMomentumMap}
    J_\xi \star J_\eta - J_\eta \star J_\xi 
    = \I\lambda J_{[\xi, \eta]} 
\end{equation}
for all $\xi, \eta \in \mathfrak{g}$. The we consider the `rescaled'
universal enveloping algebra $U(\mathfrak{g}_{\I\lambda})$ where we
rescale the Lie bracket of $\mathfrak{g}$ by $\I\lambda$ and view it
as Lie algebra over $\mathbb{C}[[\lambda]]$. Then $J$ extends to a
momentum map for the Lie algebra action induced by the
$\star$-commutator with $J_\xi$. This deforms the classical action in
the following sense: the lowest order terms in the star product
commutator are
\begin{equation}
    \label{eq:ClassLimAction}
    [J_\xi, f]_\star 
    = \I \lambda \{J_\xi, f\} + \cdot
    = -\I\lambda \xi_M f + \cdots,
\end{equation}
where $\xi_M$ is the fundamental vector field on $M$ corresponding to
$\xi$. Clearly, the action $f \mapsto [J_\xi, f]_\star$ is inner by
the very construction. Note however, that the classical action $f
\mapsto \xi_M f$ is not at all inner (unless it is trivial).

Star products with such quantum momentum maps have been excessively
studied, see e.g.\cite{gutt.rawnsley:2003a,
  mueller-bahns.neumaier:2004a}. However, if we pass to a Lie group
action with a Lie group acting on the underlying manifold by
diffeomorphisms, such an action can never be inner neither classically
nor for a star product. This is easy to see as the right hand side of
\eqref{eq:MomentumMap} for star products is local in the algebra
element $a$ while the action of a diffeomorphism via pull-backs
certainly not.

%
%

\section{The representation theories}
\label{sec:RepresentationTheories}

As we consider algebras $\mathcal{A}$ with quite a lot additional
structures the purely algebraic framework of modules will not be the
appropriate category for interesting representations. Thus we look for
representation spaces with extra structure reflecting the presence of
a $^*$-involution and symmetry of the algebras. The following
definition turns out to be very useful:
\begin{definition}[Inner Product Module]
    \label{definition:InnerProductModule}
    An \emph{inner product module} $\EA$ over $\mathcal{A}$ is a right
    $\mathcal{A}$-module with inner product
    \[
    \SPA{\cdot, \cdot}: \EA \times \EA \longrightarrow \mathcal{A}
    \]
    such that for all $x, y \in \EA$ and $a \in \mathcal{A}$
    \begin{enumerate}
    \item $\SPA{\cdot, \cdot}$ is $\ring{C}$-linear in second argument,
    \item $\SPA{x, y} = {\SPA{y, x}}^*$,
    \item $\SPA{x, y \cdot a} = \SP{x, y} a$,
    \item $\SPA{\cdot, \cdot}$ is non-degenerate.
    \end{enumerate}
    In addition, $\EA$ is called \emph{$H$-covariant} if it carries a
    $H$-action with
    \begin{equation}
        \label{eq:gActsOnSPA}
        g \acts \SPA{x, y} 
        = \SPA{S(g_\sweedler{1})^* \acts x, g_\sweedler{2} \acts y}.
    \end{equation}
\end{definition}
It follows from the non-degeneracy that $g \acts (x \cdot a) =
(g_\sweedler{1} \acts x) \cdot (g_\sweedler{2} \acts a)$.

Having an inner product module one can consider the
following operators on it. We call a map $A : \EA \longrightarrow \EpA$
from one inner product right $\mathcal{A}$-module to another
\emph{adjointable} if there exists a map $A^*: \EpA \longrightarrow
\EA$ with
\begin{equation}
    \label{eq:Adjointable}
    \SPEpA{x', Ax} = \SPEA{A^*x', x}
\end{equation}
for all $x \in \EA$ and $x' \in \EpA$. It is then an easy check that
the maps $A$ and $A^*$ are necessarily right $\mathcal{A}$-linear,
$A^*$ is uniquely determined by $A$ and $A^*$ is adjointable as well
with $A^{**} = A$. Denoting the set of all adjointable operators by
$\mathfrak{B}(\EA, \EpA)$ we have the usual properties: linear
combinations and composition of adjointable operators are again
adjointable with adjoints given in the usual way. In particular,
$\mathfrak{B}(\EA) = \mathfrak{B}(\EA, \EA)$ becomes a unital
$^*$-algebra itself. If $\EA$ is $H$-covariant then
$\mathfrak{B}(\EA)$ carries a $^*$-action of $H$ itself, given
explicitly by
\begin{equation}
    \label{eq:InducedStarActionOfBEA}
    (g \acts A)x 
    = g_\sweedler{1} \acts (A S(g_\sweedler{2}) \acts x)
\end{equation}
for $x \in \EA$ and $A \in \mathfrak{B}(\EA)$. Note however that this
is \emph{not} an inner action. The reason is that the maps $x \mapsto
g \acts x$ are \emph{not} adjointable in general and hence not given
by elements of $\mathfrak{B}(\EA)$. This follows immediately from
\eqref{eq:gActsOnSPA}: only if the action $g \acts \SPEA{x, y}$ would
be the trivial action, the map $x \mapsto g \acts x$ is adjointable
with adjoint given by $x \mapsto g^* \acts x$.

Having an inner product right $\mathcal{A}$-module $\EA$ we can define
$^*$-representations of a $^*$-algebra $\mathcal{B}$ on $\EA$ as
follows:
\begin{definition}[$^*$-Representation]
    \label{definition:Rep}
    A \emph{$^*$-representation} of $\mathcal{B}$ on $\EA$ is a
    $^*$-homomorphism $\mathcal{B} \longrightarrow \mathfrak{B}(\EA)$.
    It is called $H$-covariant if in addition
    \[
    g \acts (b \cdot x) = (g_{(1)} \acts b) \cdot (g_{(2)} \acts x)
    \]
    An \emph{intertwiner} is an adjointable $\mathcal{B}$-module map
    $T: \EA \longrightarrow \EpA$. A \emph{$H$-covariant intertwiner}
    is an intertwiner which is $H$-equivariant.
\end{definition}
In other words, a $H$-covariant $^*$-representation is a
$^*$-homomorphism which intertwines the $^*$-action of $H$ on
$\mathcal{B}$ with the canonical $^*$-action on $\mathfrak{B}(\EA)$.
\begin{definition}[$^*$-Representation Theory]
    \label{definition:sModCategory}
    Let $\mathcal{A}$, $\mathcal{B}$ be $^*$-algebras over $\ring{C}$.
    \begin{itemize}
    \item The category of $^*$-representations of $\mathcal{B}$ on
        inner product right $\mathcal{A}$-modules with intertwiners as
        morphisms is denoted by $\sMod[\mathcal{A}](\mathcal{B})$.
    \item The category of $H$-covariant $^*$-representations of
      $\mathcal{B}$ on inner product right $\mathcal{A}$-modules with
      $H$-covariant intertwiners as morphisms is denoted by $\sMod[H,
      \mathcal{A}](\mathcal{B})$.
  \end{itemize}
\end{definition}
It is clear that we indeed obtain categories this way as the
composition of intertwiners is again an intertwiner in both
situations.
\begin{remark}
    \label{remark:Unital}
    To be more precise, we should add that the above definitions
    require \emph{unital} $^*$-algebras and modules where the algebra
    unit always acts as identity. For reasons to become clear later
    on, the case of non-unital $^*$-algebras is much more technical
    and one should add the condition that $\mathcal{B} \cdot \BEA =
    \BEA$ in the above definition. However, we shall mainly discuss
    unital algebras in the sequel.
\end{remark}

Up to now, we have not yet used the positivity available from the
ordering of $\ring{R}$ at all. To incorporate this physically most
important feature one first defines positive functionals as follows: a
linear functional $\omega: \mathcal{A} \longrightarrow \ring{C}$ is
called \emph{positive} if
\begin{equation}
    \label{eq:omegaPositive}
    \omega(a^*a) \ge 0
\end{equation}
for all $a \in \mathcal{A}$. It is easy to see that $\omega$ satisfies
a Cauchy-Schwarz inequality and the reality condition $\omega(a^*b) =
\cc{\omega(b^*a)}$ which gives $\omega(a^*) = \cc{\omega(a)}$ if
$\mathcal{A}$ is e.g. unital. Having positive functionals one defines
\emph{positive algebra elements} to be those elements $a \in
\mathcal{A}$ with $\omega(a) \ge 0$ for all positive functionals
$\omega$. They form a convex cone in $\mathcal{A}$ stable under the
operations $a \mapsto b^* a b$ for arbitrary $b \in \mathcal{A}$.  The
cone of positive elements is denoted by $\mathcal{A}^+$ while the
convex cone elements of the form $\sum_i \lambda_i a_i^*a_i$ with
$\lambda_i > 0$ and $a_i \in \mathcal{A}$ is denote by
$\mathcal{A}^{++}$. We have $\mathcal{A}^{++} \subseteq \mathcal{A}^+$
and in general the inclusion is strict. A remarkable exception are the
$C^*$-algebras where $\mathcal{A}^{++} = \mathcal{A}^+$ and in fact
every positive $a$ can uniquely be written as $a = b^2$ with a
positive $b$, the square root of $a$.  More sophisticated notions of
positivity are discussed in \cite{schmuedgen:1990a} based on
particular choices of sub-cones of positive functionals, see
\cite{waldmann:2004a} for a comparison. However, here we only use the
above definition.

Analogously to the situation of Hilbert modules over $C^*$-algebras
one defines in our completely algebraic framework pre-Hilbert modules:
\begin{definition}[Pre-Hilbert module]
    \label{definition:PreHilbertModule}
    An inner product right $\mathcal{A}$-module $\EA$ is called
    \emph{pre-Hilbert module} if for all $n$ and $x_1, \ldots, x_n \in
    \EA$
    \begin{equation}
        \label{eq:CompletelyPositive}
        \left(\SPA{x_i, x_j}\right) \in M_n(\mathcal{A})^+
    \end{equation}
    i.e. $\SPA{\cdot, \cdot}$ is \emph{completely positive}.
\end{definition}
In fact, for a $C^*$-algebra, a positive inner product is always
completely positive. In general, we have to use the above definition
in order to have good behaviour under tensor products. Restricting to
pre-Hilbert modules instead of general inner product modules gives
more specific notions for $^*$-representations:
\begin{itemize}
\item The sub-category of $^*$-representations of $\mathcal{B}$ on
    pre-Hilbert modules over $\mathcal{A}$ is denoted by
    $\Rep[\mathcal{A}](\mathcal{B})$.
\item The sub-category of $H$-covariant $^*$-representations of
    $\mathcal{B}$ on pre-Hilbert modules over $\mathcal{A}$ is
    denoted by $\Rep[H, \mathcal{A}](\mathcal{B})$.
\end{itemize}
The most important case will be when the auxiliary $^*$-algebra
$\mathcal{A}$ is just the ring of scalars $\ring{C}$ (with trivial
action of $H$).

In principle, for a given $^*$-algebra $\mathcal{B}$, one would like
to understand the category $\sMod[\mathcal{A}](\mathcal{B})$ for any
coefficient algebra $\mathcal{A}$ or at least for $\mathcal{A} =
\ring{C}$ is some detail: basic questions would be to understand the
``irreducible'' representations, a notion for which it is not even
clear what is appropriate, the decomposition of a given representation
into irreducible ones, etc. Of course, in this generality and even for
$^*$-algebras like $\mathcal{B} = (C^\infty(M)[[\lambda]], \star)$
such a program is much to hard to be attacked successfully. Instead,
one has to be more modest and try to find ``interesting''
representations, e.g. by a GNS construction out of ``interesting''
positive functionals. In the case of star products ``interesting''
could mean that the positive functionals have some concrete geometric
and hence classical interpretation. Here many results have been found,
see e.g. the overview in \cite{waldmann:2005b}.

%
%

\section{Tensor products and Morita theory}
\label{sec:TensorProductsMoritaTheory}

In this section we shall now proceed with our investigation of the
representation theories, but from a different point of view: even
though it is hard or even impossible to describe
$\sMod[\mathcal{A}](\mathcal{B})$ or $\Rep[\mathcal{A}](\mathcal{B})$
in a reasonably ``explicit'' way, it might well be possible to
\emph{compare} the representation theories of different $^*$-algebras
and determine whether they are ``the same''.

To be more precisely, we are interested in finding functors form one
category of representations into another with particular interest in
\emph{equivalences} of categories. In the framework of rings and
modules this is the realm of Morita theory. Thus we want to adapt the
well-known ring-theoretic notions of Morita theory to our more
specific categories of modules.

The first step is to find an appropriate notion of tensor products.
Here we can rely on the ideas of Rieffel \cite{rieffel:1974a,
  rieffel:1974b} from $C^*$-algebra theory. In fact, the following
definition also makes sense in our much more algebraic framework:
Given $\CFB \in \sMod[\mathcal{B}](\mathcal{C})$ and $\BEA \in
\sMod[\mathcal{A}](\mathcal{B})$ we can define \emph{Rieffel's inner
  product} on their $\mathcal{B}$-tensor product by
\begin{equation}
    \label{eq:RieffelInnerProduct}
    \SPFEA{\phi \otimes x, \psi \otimes y}
    =
    \SPEA{x, \SPFB{\phi, \psi} \cdot y}
\end{equation}
for factorizing tensors and extend this to a $\mathcal{A}$-valued
inner product on $\mathcal{F} \tensor[\mathcal{B}] \mathcal{E}$. Since
\eqref{eq:RieffelInnerProduct} obviously has the correct linearity
properties, this extension is possible. However, the inner product may
still be degenerate. In fact, for $C^*$-algebras this can not happen,
but in general, there are examples where one has degeneracy. Thus we
have to divide by the degeneracy space
\begin{equation}
    \label{eq:TensorProduct}
    \CFB \tensM[\mathcal{B}] \BEA
    =
    (\CFB \tensor[\mathcal{B}] \BEA)
    \big/
    (\CFB \tensor[\mathcal{B}] \BEA)^\bot,
\end{equation}
and arrive at a $(\mathcal{C}, \mathcal{A})$-bimodule with
non-degenerate $\mathcal{A}$-valued inner product. Since $\mathcal{C}$
acts by adjointable operators, the left $\mathcal{C}$-module structure
survives the quotient~\eqref{eq:TensorProduct}. It is then an easy
check to verify that \eqref{eq:TensorProduct} indeed is a
$^*$-representation of $\mathcal{C}$ on an inner product right
$\mathcal{A}$-module. We call $\tensM[\mathcal{B}]$ the internal
tensor product (over $\mathcal{B}$).

As long as we do not take care of positivity we are done with the
above construction yielding a tensor product $\tensM$ with good
functorial properties, see \cite{ara:1999a}. Taking positivity into
account, it becomes more tricky: here the internal tensor product of
positive inner products is not behaving well. Instead, one needs
\emph{completely} positive inner products. Then one can show that
their tensor product \eqref{eq:RieffelInnerProduct} gives again a
completely positive inner product \cite{bursztyn.waldmann:2001a,
  bursztyn.waldmann:2005b}.

Finally, we have to consider the $H$-covariant situation: this is
again easy, one just computes that if both inner products on $\CFB$
and $\BEA$ satisfy \eqref{eq:gActsOnSPA} then the Rieffel inner
product is $H$-covariant again with respect to the canonical
$H$-action on the tensor product, see \cite{jansen.waldmann:2006a,
  jansen:2006a}. Collecting all these results together with the afore
mentioned functoriality properties eventually gives the following
theorem \cite{ara:1999a, bursztyn.waldmann:2005b,
  bursztyn.waldmann:2001a, jansen.waldmann:2006a, jansen:2006a}:
\begin{theorem}[Internal Tensor Product]
    \label{theorem:InternalTensorProduct}
    The internal tensor product $\tensM$ gives functors
    \begin{equation}
        \label{eq:TensorsMod}
        \tensM: 
        \sMod[\mathcal{B}](\mathcal{C}) \times 
        \sMod[\mathcal{A}](\mathcal{B}) \longrightarrow
        \sMod[\mathcal{A}](\mathcal{C})
    \end{equation}    
    \begin{equation}
        \label{eq:TensorRep}
        \tensM: 
        \Rep[\mathcal{B}](\mathcal{C}) \times 
        \Rep[\mathcal{A}](\mathcal{B}) \longrightarrow
        \Rep[\mathcal{A}](\mathcal{C})
    \end{equation}
    \begin{equation}
        \label{eq:TensorsModH}
        \tensM: 
        \sMod[H, \mathcal{B}](\mathcal{C}) \times 
        \sMod[H, \mathcal{A}](\mathcal{B}) \longrightarrow
        \sMod[H, \mathcal{A}](\mathcal{C})
    \end{equation}
    \begin{equation}
        \label{eq:TensorRepH}
        \tensM: 
        \Rep[H, \mathcal{B}](\mathcal{C}) \times 
        \Rep[H, \mathcal{A}](\mathcal{B}) \longrightarrow
        \Rep[H, \mathcal{A}](\mathcal{C})
    \end{equation}
    Moreover, $\tensM$ is associative up to the usual unitary
    intertwiners.
\end{theorem}

The theorem can now be used to enhance the category of $^*$-algebras
in the following way: instead of taking $^*$-homomorphisms as
morphisms between $^*$-algebras we take isometric isomorphism classes
of $^*$-representations $\BEA \in \sMod[\mathcal{A}](\mathcal{B})$ as
morphisms from $\mathcal{A}$ to $\mathcal{B}$. The composition of such
bimodules is then the internal tensor product $\tensM$. Then unit
arrow will be the canonical $^*$-representation $\AAA$ with inner
product $\SP{a, b} = a^*b$. Since the tensor product $\tensM$ is
associative up to unitary intertwiners we get indeed an associative
composition law on the level of isomorphism classes. The only
difficulty is that the bimodule $\AAA$ might not act as unit at all.
The problem is that if $\mathcal{A}$ does not have a unit element,
then tensoring with $\AAA$ might not give a result isomorphic to the
representation we started with. Examples are easily found. The way out
is to restrict either to unital $^*$-algebras, or, more generally, to
$^*$-algebras which are non-degenerate and idempotent: non-degenerate
means that $ab = 0$ for all $a$ implies $b = 0$ and idempotent means
that products $ab$ span already $\mathcal{A}$. In the following, we
shall stick to unital $^*$-algebras for convenience, the details for
the non-degenerate and idempotent case can be found in
\cite{bursztyn.waldmann:2005b}.

With these considerations, the following definition makes sense: for
$\mathcal{A}$, $\mathcal{B}$ we define the class
\begin{equation}
    \label{eq:BiMods}
    \BiMods(\mathcal{B}, \mathcal{A})
    =
    \left\{
        \textrm{isomorphism classes of} \; 
        \BEA \in \sMod[\mathcal{A}](\mathcal{B})
    \right\}
\end{equation}
and denote by $\BiMods$ the category whose objects are unital
$^*$-algebras over $\ring{C}$ with the isomorphism classes
$\BiMods(\mathcal{B}, \mathcal{A})$ of bimodules as morphisms from
$\mathcal{A}$ to $\mathcal{B}$, the internal tensor product $\tensM$
as composition law, and the isomorphism class of $\AAA$ as unit arrow.
It is then clear from the functoriality of $\tensM$ that $\tensM$ is
well-defined on isomorphism classes.  Therefor we end up with a
category of bimodules, completely analogous to the well-known
ring-theoretic case.

In a completely analogous fashion we can define the category
$\BiModsH$ of $^*$-algebras with isomorphism classes of $H$-covariant
bimodules as morphisms, the category $\BiModstr$ of $^*$-algebras with
isomorphism classes of bimodules with completely positive inner
products as morphisms, and finally, the category $\BiModstrH$ of
$^*$-algebras with isomorphism classes of $H$-covariant bimodules with
completely positive inner products as morphisms.

\begin{remark}
    \label{remark:BicategoryApproach}
    In fact, identifying bimodules up to isomorphisms might not be a
    good idea after all. Instead, one can equally well take \emph{all}
    bimodules (of the particular, specific types). Then the tensor
    product $\tensM$ fails to be associative and $\AAA$ fails to be a
    unit arrow. However, the failure is encoded in a functorial way
    whence we end up with a \emph{bicategory} in each of the above
    four cases. The $1$-morphisms are the bimodules and the
    $2$-morphisms, i.e. the arrows between the arrows, are the
    intertwiners, always in the correct version respecting the inner
    products and the $H$-covariance. A detailed study of these
    bicategorical aspects can be found in \cite{waldmann:2008a:script}
    as well as in \cite{jansen:2006a}. Its ring-theoretic version goes
    back to Benabou \cite{benabou:1967a} and was probably the
    motivation to study bicategories at all.
\end{remark}

Using either the bicategorical approach or the one presented above,
Morita equivalence is now simply isomorphism in the enhanced
categories: 
\begin{definition}[Morita equivalence]
    \label{definition:MoritaEquvialence}
    Two (unital) $^*$-algebras $\mathcal{A}$ and $\mathcal{B}$ are called
    \begin{itemize}
    \item \emph{$^*$-Morita equivalent} if they are isomorphic in
        $\BiMods$.
    \item \emph{Strongly Morita equivalent} if they are isomorphic in
        $\BiModstr$.
    \item \emph{$H$-covariantly $^*$-Morita equivalent} if they are
        isomorphic in $\BiModsH$.
    \item \emph{$H$-covariantly strongly Morita equivalent} if they
        are isomorphic in $\BiModstrH$.
    \end{itemize}
    A bimodule representing such isomorphism is called
    \emph{equivalence bimodule}.
\end{definition}
The case of $C^*$-algebras was first studied by Rieffel in
\cite{rieffel:1974b} and serves as motivation for all the
generalizations. In fact, he coined the name strong Morita
equivalence. Ara discussed the notion of $^*$-Morita equivalence in
\cite{ara:1999a} while in \cite{bursztyn.waldmann:2001a,
  bursztyn.waldmann:2005b} we extended Rieffel's notion to arbitrary
$^*$-algebras now taking into account the complete positivity of inner
product compared to the approach of Ara. Finally, the $H$-covariant
situation was discussed in detail in \cite{jansen.waldmann:2006a,
  jansen:2006a} based on earlier formulations in $C^*$-algebra theory
where the symmetry was implemented by a strongly continuous group
action of some suitable topological group.

While the above definition of Morita equivalence is very elegant it
needs further work to get a formulation suitable for practical
purposes. In fact, one has to find criteria whether a given bimodule
is ``invertible'' or not. In the ring-theoretic version the
equivalence bimodules are the full projective modules over
$\mathcal{A}$ such that $\mathcal{B}$ is isomorphic to the right
$\mathcal{A}$-linear endomorphisms. However, in this classical
formulation the fact that one deals with \emph{unital} rings is
crucial. The following formulation based on Rieffel's pioneering work
avoids the usage of the unit and extends to non-degenerate and
idempotent $^*$-algebras as well. The $^*$-Morita version is due to
Ara, the strong version is from \cite{bursztyn.waldmann:2005b}:
\begin{theorem}
    \label{theorem:MoritaEquivalence}
    The following statements are equivalent:
    \begin{itemize}
    \item $\BEA \in \BiMods(\mathcal{B}, \mathcal{A})$ is an
        equivalence bimodule.
    \item On $\BEA \in \sMod[\mathcal{A}](\mathcal{B})$ there is a
        $\mathcal{B}$-valued inner product $\BSPE{\cdot, \cdot}$ such
        that
        \begin{enumerate}
        \item $\SPEA{\cdot, \cdot}$ is full, i.e. $\SPEA{\mathcal{E},
              \mathcal{E}} = \mathcal{A}$.
        \item $\BSPE{\cdot, \cdot}$ is full.
        \item $\BSPE{x, y} \cdot z = x \cdot \SPEA{y, z}$ for all $x,
            y, z \in \BEA$.
        \end{enumerate}
    \end{itemize}
    The analogous statement holds for strong Morita equivalence. Here,
    in addition, $\BSPE{\cdot, \cdot}$ is completely positive.
\end{theorem}
The unital version of this theorem is a rather straightforward
adaption of the ring-theoretic Morita theorem. To include
non-degenerate and idempotent $^*$-algebras one has to put in some
more effort.

Not surprising, there is also a characterization of the $H$-covariant
equivalence bimodules, both for the $^*$-equivalence and the strong
equivalence bimodules \cite{jansen.waldmann:2006a}:
\begin{theorem}
    \label{theorem:HcovariantMorita}
    The following statements are equivalent:
    \begin{itemize}
    \item $\BEA \in \BiModsH(\mathcal{B}, \mathcal{A})$ is an
        equivalence bimodule.
    \item On $\BEA \in \sMod[H, \mathcal{A}](\mathcal{B})$ there is a
        $\mathcal{B}$-valued inner product $\BSPE{\cdot, \cdot}$ such
        that
        \begin{enumerate}
        \item $\BSPE{\cdot, \cdot}$ is compatible with $H$-action.
        \item $\SPEA{\cdot, \cdot}$ is full, ie. $\SPEA{\mathcal{E},
              \mathcal{E}} = \mathcal{A}$.
        \item $\BSPE{\cdot, \cdot}$ is full.
        \item $\BSPE{x, y} \cdot z = x \cdot \SPEA{y, z}$.
        \end{enumerate}
    \end{itemize}
    Again, there is an analogous statement for $H$-covariantly strong
    Morita equivalence. Here $\BSPE{\cdot, \cdot}$ is in addition
    completely positive.
\end{theorem}
Having Theorem~\ref{theorem:MoritaEquivalence}, the results in
Theorem~\ref{theorem:HcovariantMorita} are obtained rather easily by
observing that the $H$-action $\BEA$ is necessarily compatible with
the inner product $\BSPE{\cdot, \cdot}$ thanks to the compatibility of
the two inner products $\BSPE{x, y} \cdot z = x \cdot \SPEA{y, z}$.

\begin{remark}[Finite rank operators]
    \label{remark:FiniteRank}
    For a general inner product right $\mathcal{A}$-module $\EA$ one
    defines the \emph{rank-one operator} $\Theta_{x, y}: \EA
    \longrightarrow \EA$ by
    \begin{equation}
        \label{eq:Thetaxy}
        \Theta_{x, y}(z) = x \cdot \SPEA{y, z},
    \end{equation}
    where $x, y, z \in \EA$. Then it is easy to see that $\Theta_{x,
      y}$ is adjointable with $\Theta_{x, y}^* = \Theta_{y, x}$.  The
    linear span of all rank-one operators is denoted by
    $\mathfrak{F}(\EA)$ and called the \emph{finite rank operators} on
    $\EA$.  Clearly, $\mathfrak{F}(\EA) \subseteq \mathfrak{B}(\EA)$
    is a $^*$-ideal. One can now show that for a Morita equivalence
    bimodule $\BEA$ we have $\mathcal{B} \cong \mathfrak{F}(\EA)$ via
    the left action. Moreover, the $\mathcal{B}$-valued inner product
    $\BSPE{\cdot, \cdot}$ corresponds to $\Theta_{\cdot, \cdot}$ under
    this $^*$-isomorphism. Finally, in the $H$-covariant situation,
    the $H$-action on $\mathcal{B}$ corresponds to the canonical
    action \eqref{eq:InducedStarActionOfBEA} which is easily shown to
    preserve $\mathfrak{F}(\EA)$. In the unital case we have
    $\mathcal{B} \cong \mathfrak{F}(\EA) = \mathfrak{B}(\EA)$.
\end{remark}

%
%

\section{Picard groupoids}
\label{sec:PicardGroupoids}

Having realized Morita equivalence as the notion of isomorphism in the
enhanced categories $\BiMods$, $\BiModstr$, $\BiModsH$, and
$\BiModstrH$, respectively, it is of course not only interesting to
ask \emph{whether} two $^*$-algebras are Morita equivalent, but also
\emph{in who many ways}. Answers to both questions are encoded in the
Picard groupoid, for which we again have four different flavors:
\begin{definition}[Picard groupoids]
    \label{definition:PicardGroupoid}
    \begin{itemize}
    \item The (large) groupoid of invertible arrows in $\BiMods$ is
        called the \emph{$^*$-Picard groupoid} $\starPic$. The
        isotropy group at $\mathcal{A}$ is called the
        \emph{$^*$-Picard group} $\starPic(\mathcal{A})$ of
        $\mathcal{A}$.
    \item Analogously, one defines the \emph{strong Picard groupoid}
        $\StrPic$, the \emph{$H$-covariant $^*$-Picard groupoid}
        $\starPicH$, and the \emph{$H$-covariant strong Picard
          groupoid} $\StrPicH$ together with the Picard groups
        $\StrPic(\mathcal{A})$, $\starPicH(\mathcal{A})$, and
        $\StrPicH(\mathcal{A})$, respectively.
    \end{itemize}
\end{definition}
From this point of view, Morita theory consists in the following two
principle tasks:
\begin{enumerate}
\item Determine the \emph{orbits} of the Picard groupoid. These are
    precisely the Morita equivalence classes of $^*$-algebras.
\item Determine the \emph{isotropy groups} of the Picard groupoid.
    They encode in how many ways a $^*$-algebra can be Morita
    equivalent to itself.
\end{enumerate}
In particular, by groupoid abstract non-sense isotropy groups are
always isomorphic along the orbits whence the Picard groups are
isomorphic for Morita equivalent $^*$-algebras. Needless to say, the
above program should be carried through in all four versions of Morita
theory. Moreover, as unital $^*$-algebras are particular types of
rings, we also have the ring-theoretic version of Morita theory at
hand, as well as a ring-theoretic version of $H$-covariant Morita
theory. Thus we also can study the Picard groupoids $\Pic$ and
$\PicH$, respectively.

One strategy to learn something about the various Picard group(oids)
in general is to use the following groupoid morphisms which simply
forget the additional structure successively. This gives a commuting
diagram of groupoid morphisms
\begin{equation}
    \label{eq:CommutingGroupoidMorphisms}
      \bfig
      \Vtriangle(0,350)<500,200>[\StrPicH`\starPicH`\Pic_H;``]
      \Vtriangle(0,0)/@{>}|\hole`>`>/<500,200>[\StrPic`\starPic,`\Pic;``]
      \morphism(0,525)<0,-250>[`;]
      \morphism(500,325)<0,-250>[`;]
      \morphism(1000,525)<0,-250>[`;]
      \efig
\end{equation}
for each of which one would like to know kernel and image. In
\cite{bursztyn.waldmann:2005b} the arrow $\StrPic \longrightarrow
\Pic$ was studied in detail. The arrows from and to the $^*$-versions
encode how many inner products with \emph{different signature} than
the completely positive one can have. From that point of view, they
are the simplest to understand.

In the following, we concentrate on the arrow
\begin{equation}
    \label{eq:StrPicHtoStrPic}
    \StrPicH \longrightarrow \StrPic,
\end{equation}
following \cite{jansen.waldmann:2006a}. To this end, we restrict to
unital $^*$-algebras in the sequel. The non-unital case is much more
mysterious.  In order to determine the kernel of this groupoid
morphism we need some preparation. First recall that the linear maps
$\Hom_{\ring{C}}(H, \mathcal{A})$ form an associative algebra over
$\ring{C}$ with respect to the \emph{convolution} product
\begin{equation}
    \label{eq:Convolution}
    (\twist{a} * \twist{b})(g) 
    = \twist{a}(g_\sweedler{1}) \twist{b}(g_\sweedler{2})
\end{equation}
for $\twist{a}, \twist{b} \in \Hom_{\ring{C}}(H, \mathcal{A})$ and
with unit given by $\twist{e}(g) = \epsilon(g) \Unit_{\mathcal{A}}$.
We now consider the following particular linear maps:
\begin{definition}
    \label{definition:FunnyGroups}
    Define $\mathrm{GL}(H, \mathcal{A}) \subseteq \Hom_{\ring{C}}(H,
    \mathcal{A})$ to be the subset of those $\twist{a}$ with
    \begin{itemize}
    \item $\twist{a}(\Unit_H) = \Unit_\mathcal{A}$,
    \item $\twist{a}(gh) = \twist{a}(g_\sweedler{1}) (g_\sweedler{2}
        \acts \twist{a}(h))$,
    \item $(g_\sweedler{1} \acts b) \twist{a}(g_\sweedler{2}) =
        \twist{a}(g_\sweedler{1}) (g_\sweedler{2} \acts b)$,
    \end{itemize}
    for all $g, h \in H$ and $b \in \mathcal{A}$.  Furthermore, define
    $\mathrm{U}(H, \mathcal{A}) \subseteq \mathrm{GL}(H, \mathcal{A})$
    by the additional condition
    \begin{itemize}
    \item $\twist{a}(g_\sweedler{1})
        \left(\twist{a}(S(g_\sweedler{2})^*)\right)^* = \epsilon(g)
        \Unit_{\mathcal{A}}$.
    \end{itemize}
\end{definition}
The following result is a rather straightforward verification:
\begin{proposition}
    \label{proposition:TheyAreReallyGroups}
    $\mathrm{GL}(H, \mathcal{A})$ is a group with respect to the
    convolution product of $\Hom_{\ring{C}}(H, \mathcal{A})$ and
    $\mathrm{U}(H, \mathcal{A})$ is a subgroup.
\end{proposition}
Let $c \in \mathrm{U}(\mathcal{Z}(\mathcal{A}))$ be unitary and
central then
\begin{equation}
    \label{eq:hatc}
    \hat{c}(g) = (g \acts c^{-1}) c
\end{equation}
defines an element $\hat{c} \in \mathrm{U}(H, \mathcal{A})$ as a
simple computation shows. It is easy to see that this gives in fact a
group homomorphism whose image is central.
\begin{proposition}
    \label{proposition:ExactSequence}
    There is an exact sequence of groups
    \begin{equation}
        \label{eq:ExactSequence}
        1 \longrightarrow
        \mathrm{U}(\mathcal{Z}(\mathcal{A}))^H
        \longrightarrow
        \mathrm{U}(\mathcal{Z}(\mathcal{A}))
        \stackrel{\widehat{\quad}}{\longrightarrow}
        \mathrm{U}(H, \mathcal{A}),
    \end{equation}
    and the image of $\mathrm{U}(\mathcal{Z}(\mathcal{A}))$ is central
    in $\mathrm{U}(H, \mathcal{A})$
\end{proposition}
Note that in general the group $\mathrm{U}(H, \mathcal{A})$ is far
from being abelian. Nevertheless, the above proposition allows to
consider the following quotient group
\begin{equation}
    \label{eq:UNullDef}
    \mathrm{U}_0(H, \mathcal{A})
    =
    \mathrm{U}(H, \mathcal{A})
    \big/ 
    \widehat{\mathrm{U}(\mathcal{Z}(\mathcal{A}))}.
\end{equation}
Note that the same construction also works without the ``unitarity''
requirement and yields a group $\mathrm{GL}_0(H, \mathcal{A}) =
\mathrm{GL}(H, \mathcal{A}) \big/
\widehat{\mathrm{GL}(\mathcal{Z}(\mathcal{A}))}$, where now we also
allow for arbitrary invertible central elements instead of unitary
ones.

In order to understand the groupoid morphism $\StrPicH \longrightarrow
\StrPic$ we assume to have a strong equivalence bimodule $\BEA$
between $\mathcal{A}$ and $\mathcal{B}$. Then the question is whether
we can find a compatible action of $H$ on $\BEA$ turning it into a
$H$-covariant equivalence bimodule. Thus the \emph{surjectivity} is a
\emph{lifting problem}. As already simple examples show, the
computation of the image will be a very difficult task. Moreover, the
result will depend strongly on the two $^*$-algebras $\mathcal{A}$ and
$\mathcal{B}$:
\begin{example}[Lifting of group actions]
    \label{example:Lifting}
    Let $\mathcal{A} = \mathcal{B} = C^\infty(M)$ be the smooth
    functions on a manifold and assume that the symmetry is given by a
    Lie algebra action of $\mathfrak{g}$. The interesting
    self-equivalence bimodules are then known to be the sections
    $\Gamma^\infty(L)$ of complex line bundles $L \longrightarrow M$,
    endowed with their usual bimodule structure. Thus the question
    whether we can endow the bimodule $\Gamma^\infty(L)$ with a
    compatible symmetry of $\mathfrak{g}$ is equivalent to the
    question whether we can lift the Lie algebra action to $L$. In
    general, this is a difficult and of course classical question in
    differential geometry whose answer is known to depend very much on
    the example.  Analogous statements are of course valid for Lie
    group actions instead of Lie algebra actions.
\end{example}

While the image is hard to determine, we can say something about the
kernel. So we assume that we have at least one lifting, i.e. there is
at least one compatible action $\acts$ of $H$ on $\BEA$. Then we have
to understand how many compatible actions we can have. The idea is to
parameterize to possible actions starting from the given $\acts$ by
$\twist{a} \in \mathrm{U}(H, \mathcal{B})$. This yields in fact a
bijection. However, it may happen that some of the actions yield
isomorphic equivalence bimodules. It turns out that they are precisely
parameterized by those $\twist{a}$ coming from
$\mathrm{U}(\mathcal{Z}(\mathcal{A}))$ via \eqref{eq:hatc}. This
eventually yields the following result, see
\cite{jansen.waldmann:2006a}:
\begin{theorem}
    \label{theorem:Kernel}
    Assume $\StrPic(\mathcal{B}, \mathcal{A})$ is non-empty. Then one
    has the alternatives:
    \begin{itemize}
    \item $\StrPicH(\mathcal{B}, \mathcal{A}) = \emptyset$.
    \item $\StrPicH(\mathcal{B}, \mathcal{A}) \longrightarrow
        \mathrm{im}(\StrPicH(\mathcal{B}, \mathcal{A}))$ is a
        principal $\mathrm{U}_0(H, \mathcal{B})$-bundle over the
        image, i.e. $\mathrm{U}_0(H, \mathcal{B})$ acts freely and
        transitively on the fibers.
    \end{itemize}
    The same holds for the $^*$-Picard groupoids. Moreover, analogous
    statements hold in the ring-theoretic setting with
    $\mathrm{GL}_0(H, \mathcal{B})$ instead of $\mathrm{U}_0(H,
    \mathcal{B})$.
\end{theorem}
By symmetry of Morita equivalence it is immediately clear that the
groups $\mathrm{U}_0(H, \mathcal{B})$ are \emph{invariant} under
$H$-covariant $^*$-equivalence while $\mathrm{GL}_0(H, \mathcal{A})$
is invariant under $H$-covariant ring-theoretic Morita invariance. In
fact, one can show that the whole exact sequence
\eqref{eq:ExactSequence} is a \emph{Morita invariant}.

Back in our geometric context, one can actually compute the group
$\mathrm{U}_0(H, \mathcal{A})$ in terms of classical geometric data:
\begin{example}
    \label{example:UNullGeometric}
    Let $\mathcal{A} = C^\infty(M)$ and let $H =
    U_{\mathbb{C}}(\mathfrak{g})$ be the complexified universal
    enveloping algebra of a real Lie algebra $\mathfrak{g}$ acting on
    $M$ by vector fields. Then one can show \cite{waldmann:2006a}
    \begin{equation}
        \label{eq:UNullGeometric}
        \mathrm{U}_0(H, \mathcal{A}) = 
        H^1_{CE}(\mathfrak{g}, C^\infty(M, \I \mathbb{R})) 
        \big/
        \widehat{H^1_{dR}(M, 2\pi\I \mathbb{Z})},
    \end{equation}
    where $H^1_{dR}(M, 2\pi\I \mathbb{Z})$ denotes the
    $2\pi\I$-integral deRham cohomology of $M$ and the map
    \begin{equation}
        \label{eq:LogDerivative}
        \widehat{\quad} : H^1_{dR}(M, 2\pi\I \mathbb{Z}) \longrightarrow
        H^1_{CE}(\mathfrak{g}, C^\infty(M, \I\mathbb{R}))
    \end{equation}
    is obtained from a ``logarithmic derivative'' induced by
    \eqref{eq:hatc}. Finally, $H^1_{CE}(\mathfrak{g}, C^\infty(M, \I
    \mathbb{R}))$ denotes the Chevalley-Eilenberg cohomology of
    $\mathfrak{g}$ with values in the imaginary-valued functions on
    $M$. Thus we obtain a linear space divided by some integer lattice
    for the group $\mathrm{U}_0(H, \mathcal{A})$.
\end{example}
\begin{remark}
    \label{remark:OtherAlgebras}
    We note that the above result generalizes to other algebras
    $\mathcal{A}$ which have a sort of ``exponential function'', see
    \cite{waldmann:2006a} for details. We also note that the above
    statement gives a classification of how many inequivalent lifts of
    the Lie algebra action on $M$ one has for a complex line bundle.
    This way, one can re-produce well-known results from differential
    geometry (the lifting problem is of course classical) by Morita
    theoretic considerations.
\end{remark}

%
%

\section{Further remarks}
\label{sec:FurtherRemarks}

We have seen that the surjectivity of \eqref{eq:StrPicHtoStrPic} is
quite hard to understand in general. However, under the additional
assumption of \emph{inner actions} of $H$ on $\mathcal{A}$ and
$\mathcal{B}$ one can always lift: indeed, assume $\BEA$ is a strong
Morita equivalence bimodule then
\begin{equation}
    \label{eq:TheLiftForInner}
    g \acts x 
    =
    J_{\mathcal{B}}(g_{(1)})\cdot x \cdot J_{\mathcal{A}}(S(g_{(2)}))
\end{equation}
defines a compatible action of $H$ on $\BEA$. In fact, one even has
the following statement, see the forthcoming project
\cite{jansen.neumaier.waldmann:2005a:pre}:
\begin{theorem}
    \label{theorem:LiftIsCool}
    Assume $\StrPic(\mathcal{B}, \mathcal{A}) \ne \emptyset$ and let
    both algebras have inner actions.  Then 
    \begin{equation}
        \label{eq:StrPicHSurjective}
        \StrPicH(\mathcal{B}, \mathcal{A})
        \longrightarrow \StrPic(\mathcal{B}, \mathcal{A})
    \end{equation}
    is surjective. Moreover,
    \[
    1 \longrightarrow \mathrm{U}_0(H, \mathcal{A})
    \longrightarrow \StrPicH (\mathcal{A})
    \longrightarrow \StrPic(\mathcal{A})
    \longrightarrow 1
    \]
    splits via the momentum map.
\end{theorem}
Thus we are in some sense in an optimal situation: we have determined
the $H$-covariant Picard groupoid completely in terms of the groups
\eqref{eq:UNullDef} and the usual Picard groupoid. Of course, having
inner actions is a severe restriction in general. Nevertheless, there
are important examples as discussed in Section~\ref{sec:SetUp}.

In \cite{jansen.neumaier.waldmann:2005a:pre}, we investigate the
particular case of \emph{symplectic} star product algebras where the
existence of a quantum momentum map is well-understood. In particular,
one can use these results to obtain lifting results this way.

At last, we remark that the study of the Picard groupoids for star
product algebras is understood best if one uses the classical limit
functor, where the formal parameter $\lambda$ is set to zero. Then one
can study how the Picard groupoid behaves under formal deformations of
the underlying algebras. As shown in \cite{bursztyn.waldmann:2004a}
this allows to explicitly compute the ring-theoretic Picard groups of
certain star product algebras as well as the orbits
\cite{bursztyn.waldmann:2002a}. Thus it remains as a challenge to use
these results and extend them to the $H$-covariant as well as to the
$^*$-equivalence and strong equivalence setting.

%
%


\end{document}